\documentclass[11pt, draft]{amsart}
\usepackage{amssymb}

\theoremstyle{plain}
\newtheorem{theorem}{Theorem}

\theoremstyle{remark}
\newtheorem{remark}{Remark}

\begin{document}
\title[]{An extension of a characterization of the automorphisms of
Hilbert space effect algebras}
\author{LAJOS MOLN\'AR}
\address{Institute of Mathematics and Informatics\\
         University of Debrecen\\
         4010 Debrecen, P.O. Box 12, Hungary}
\email{molnarl@math.klte.hu}
\author{ENDRE KOV\'ACS}
\address{Department of Theoretical Physics\\
         University of Debrecen\\
         4010 Debrecen, P.O. Box 5, Hungary}
\email{kendre@dtp.atomki.hu}
\thanks{This paper was written when the first author held a Humboldt
Research Fellowship. He is very grateful to
the Alexander von Humboldt Foundation for providing ideal conditions
for research and also to his host Werner Timmermann (TU Dresden, Germany)
for very warm hospitality.
         The first author also acknowledges support from
         the Hungarian National Foundation for Scientific Research
         (OTKA), Grant No. T030082, T031995, and from
         the Ministry of Education, Hungary, Grant
         No. FKFP 0349/2000}
\date{October 27, 2002}
\begin{abstract}
The aim of this paper is to
show that if an order preserving bijective transformation of the
Hilbert space effect algebra also preserves the probability with
respect to a fixed pair of mixed states, then it is an ortho-order
automorphism. A similar result for the orthomodular lattice of all
sharp effects (i.e., projections) is also presented.

\end{abstract}
\maketitle

Effect algebras play fundamental role in the theory of quantum
measurement \cite{BusLahMit} (also see \cite{Kraus} and \cite{Ludwig}).
In the Hilbert space setting, the so-called Hilbert space
effect algebra ${\mathcal E}(H)$ is the set of all positive
bounded linear operators
on the Hilbert space $H$ which are majorized by the identity $I$. This
set is usually equipped with certain operations and/or relations which
all have physical meaning. Therefore, there are different algebraic
structures on ${\mathcal E}(H)$.
In some respect, probably the most important
such structure is obtained when we equip ${\mathcal E}(H)$ with the
partial order $\leq$ (which is just the usual order among self-adjoint
operators restricted to ${\mathcal E}(H)$) and a kind of
orthocomplementation, namely, $\perp: E\mapsto
I-E$.

The study of the automorphisms of given algebraic structures is
a very important general problem in mathematics. As for
effect algebras,
the investigation of the so-called ortho-order automorphisms of
${\mathcal E}(H)$ (that is, the automorphisms with respect to the order
and orthocomplementation) was begun by Ludwig in
\cite{Ludwig}. In fact, in \cite[Section V.5]{Ludwig}
he showed that, in case $\dim H\geq 3$, these automorphisms are
implemented
by unitary-antiunitary operators. However, his argument seemed
to contain some gaps and the proof
was recently clarified in \cite{CVLL2} completely. (We mention that in
our recent
paper \cite{MoPa} we have shown that Ludwig's result holds also in the
2-dimensional case and this answers a question that was open for quite
some time.)

In our paper \cite{Mo2} we initiated the study of the automorphisms of
effect algebras (or any other quantum structure) by means of their
preserver properties. We expressed our belief there that, similarly
to the case of linear preserver problems in matrix theory (concerning
which we refer, for
example, to the survey papers \cite{LiTs, LiPi}), such investigations
may give important new
information about the automorphisms in question and they may help to
better understand the underlying algebraic structures.
According to this, in \cite{Mo2} we presented some characterizations of
the automorphisms of effect algebras via their preserver properties.
This study was continued in \cite{Mo3, Mo4} where
we obtained results of the same kind for
the automorphisms of the Jordan algebra of all bounded observables.

Turning to the content of the present paper, we remark that
the order is undoubtedly a very important relation on ${\mathcal E}(H)$.
One of the reasons is the following. As it turns out
from \cite{BuGu} (and, in fact, was asserted already by Ludwig), the
effects are determined by the weak atoms they majorize. As weak atoms
can be defined by the order exclusively, it is obvious how
essential the order is in the description of effects.
We next refer to \cite{MoGu} to see how "strong" this relation is.
In spite of this, the order preserving property alone is not strong
enough to
characterize the ortho-order automorphisms in even some weak sense.
In fact, if we consider the transformation
\[
E \longmapsto
\biggl( \frac{T^2}{2I-T^2} \biggr)^{-\frac{1}{2}}
\biggl( \bigl (I-T^2+T(I+E)^{-1}T\bigr )^{-1}-I\biggr )
\biggl( \frac{T^2}{2I-T^2} \biggr)^{-\frac{1}{2}}
\]
where $T\in {\mathcal E}(H)$ is fixed and invertible, one
can easily check that this is a bijective map on ${\mathcal
E}(H)$ which preserves the order but has nothing to do with the
ortho-order automorphisms of ${\mathcal E}(H)$.
(Here we mention that the situation is much different with the partially
ordered set $B_s(H)$ of all bounded observables on $H$. It might be quite
surprising that, as it turns out from our paper \cite{Mo3}, if
$\dim H>1$ and $\phi$ is an order preserving bijective map on $B_s(H)$, then
$\phi$ can be written in a nice explicit form which is closely related
to the form of the automorphisms of $B_s(H)$ as a Jordan algebra (cf.
\cite{CVLL}).
So, it is clear that beside the order, the preservation of some other
physical quantity is needed to characterize the automorphisms of
${\mathcal E}(H)$. The aim of this paper is to present such a result
which is a significant generalization of one of our former
results. Namely, in \cite[Theorem 2]{Mo2}
we proved that if $\dim H\geq 3$ and $\phi: {\mathcal E}(H)\to {\mathcal
E}(H)$ is a bijective map which preserves the order and also
preserves the probability with respect to a fixed pair of pure states,
then $\phi$ is an ortho-order automorphism. It is now a natural question
that what about the mixed states?
In which follows we show that the same conclusion holds also for
mixed states even in the case when $\dim H=2$ (the 1-dimensional
case is trivial). This generalizes and extends our result in
\cite{Mo2}.
Furthermore, we present a similar result concerning the
orthomodular lattice of projections.

As for the notation, let $H$ be a complex Hilbert space and, just as
before, let ${\mathcal E}(H)$ be the set of all self-adjoint bounded
linear operators $E$ on $H$ for which $0\leq E\leq I$. Denote by
${\mathcal
P}(H)$ the set of all projections on $H$. The elements of ${\mathcal
E}(H)$ are called effects (or, in other terminology, unsharp events),
while the elements of ${\mathcal P}(H)$ are called sharp effects (or, in
other terminology, sharp events). By a (mixed) state we mean a
positive trace-class operator on $H$ with trace 1. The usual
trace-functional is denoted by $\text{tr}$. The states form
a convex set whose extreme points are exactly the rank-one projections
which are called pure states.
Finally, we emphasize that when we say that a transformation preserves a
certain relation, we always mean that this relation is preserved in both
directions.

Now, the main result of the paper reads as follows.

\begin{theorem}\label{T:kov2}
Let $\phi : {\mathcal E}(H)\to {\mathcal E}(H)$ be a bijective map with the property that
\[
E \leq F \Longleftrightarrow \phi(E)\leq\phi(F) \qquad (E,F \in
{\mathcal E}(H))
\]
and suppose that there are states $D$, $D'$ such that
\[
\text{tr}\,(\phi(E)D')=\text{tr}\, (ED)
\qquad (E\in {\mathcal E}(H)).
\]
Then there exists an either unitary or antiunitary operator U on H such
that $\phi$ is of the form
\[
\phi(E)=UEU^{*} \qquad (E\in {\mathcal E}(H)).
\]
\end{theorem}

\begin{proof}
Similarly to the proof of Theorem 1 in \cite{Mo2} we obtain that,
by the order preserving property, $\phi$ preserves the projections and
their ranks.

We show that $\phi(\lambda I)=\lambda I$ and $\phi(\lambda P)=\lambda
\phi(P)$ holds for every $\lambda \in [0,1]$ and for every rank-one
projection $P$.
Let $\lambda \in [0,1]$ and pick a rank-one projection $P$ such that
$\text{tr}\, \phi(P)D'=\text{tr}\, PD\neq 0$.
As $\lambda P\leq P$, we obtain that $\phi(\lambda P)\leq \phi(P)$.
Taking into account that $\phi(P)$ is of rank 1, it follows that
$\phi(\lambda P)=\mu \phi(P)$ holds for some $\mu \in [0,1]$. Now, we
compute on one hand that
\[
\text{tr}\, \phi(\lambda P)D'=
\text{tr}\, \mu \phi(P)D'=
\mu \text{tr}\, PD
\]
and on the other hand that
\[
\text{tr}\, \phi(\lambda P)D'=
\text{tr}\, (\lambda P) D=
\lambda \text{tr}\, P D.
\]
Comparing these two equalities and remembering that $\text{tr}\, PD\neq
0$, we
deduce $\mu=\lambda$, that is, we get $\phi(\lambda P)=\lambda \phi(P)$.

At this point we need the useful concept of the strength of an effect
$E$ along a rank-one projection $Q$ (or, equivalently, along the ray
represented by the projection $Q$). According to \cite{BuGu}, this
number
is, by definition, the supremum of the set of all $t\in [0,1]$ for which
$tQ\leq E$.

Let $\lambda \in ]0,1]$. Pick a rank-one projection $P$ as above, that
is, suppose that
$\text{tr}\, \phi(P)D'=\text{tr}\, PD\neq 0$.
The strength of $\lambda I$ along $P$ is obviously $\lambda$. By the
order preserving property of $\phi$ and the homogeneity of $\phi$ on
the set of effects of the form $tP$ (see the second paragraph of our
proof), we infer that the strength of
$\phi(\lambda I)$ along $\phi(P)$ is also $\lambda$. This
holds for every rank-one projection $\phi(P)$ for which $\text{tr}\,
\phi(P)D'\neq 0$.
It is clear that we have $\text{tr}\, \phi(P)D'\neq 0$ if and only
if the range of $\phi(P)$ is not orthogonal to the range of $D'$. But
the set of all such $\phi(P)$'s is easily seen to be dense in the set of
all
rank-one projections. This implies that the strength of $\phi(\lambda
I)$ is equal to $\lambda$ along every member of a dense subset of
rank-one
projections. Lemma 5 in \cite{Mo2} tells us that in
this case we have $\phi(\lambda I)=\lambda I$.

Now, pick an arbitrary rank-one projection $P$, let $\lambda \in [0,1]$
be also arbitrary and take $\mu\in [0,1]$ such that $\phi(\lambda
P)=\mu \phi(P)$. We have on one hand that
\[
\mu \phi(P)=\phi(\lambda P)\leq \phi(\lambda I)=\lambda I
\]
implying that $\mu \leq \lambda$,
and on the other hand that
\[
\phi(\lambda P)=\mu \phi(P)\leq \mu I=\phi(\mu I)
\]
implying that $\lambda P\leq \mu I$, i.e., $\lambda \leq \mu$.
Therefore, we obtain that $\lambda =\mu$ which yields
$\phi(\lambda P)=\lambda \phi(P)$ as we have
claimed.

We next assert that $\phi$ preserves the orthogonality between rank-one
projections.
To see this, let $P,Q$ be mutually orthogonal rank-one projections.
Denote $P^\perp=I-P$
the orthogonal complement of $P$. Pick numbers $0<\lambda<\mu \leq 1$
and consider the effect $E=\lambda P+\mu P^\perp$. Clearly, we
have $\lambda
I\leq E\leq \mu I$, the strength of $E$ along $P$ is $\lambda$ and along
$Q$ (which is a subprojection of $P^\perp$) is $\mu$. It follows from
the order preserving property
of $\phi$ and from the homogeneity of $\phi$ on the scalar multiples of
rank-one projections
that $\lambda I\leq \phi(E)\leq \mu I$, the strength of $\phi(E)$ along
$\phi(P)$ is $\lambda$ and along $\phi(Q)$ is $\mu$.
Now, Lemma 3 in \cite{Mo2} tells us that in this case
the ranges of
$\phi(P)$ and $\phi(Q)$ are subspaces of the eigenspaces of $\phi(E)$
corresponding
to the eigenvalues $\lambda$ and $\mu$, respectively. This shows that
the ranges of $\phi(P)$ and $\phi(Q)$ are orthogonal to each other.

Now we have to distinguish two cases. Suppose first that $\dim H=2$. In
that case it follows from what we have just seen (where we have
made a reference to \cite[Lemma 3]{Mo2})
that
\[
\phi(\lambda P+\mu Q)=\lambda \phi(P)+\mu \phi(Q)
\]
holds for every mutually orthogonal rank-one projections $P,Q$ and real
numbers $\lambda,\mu \in [0,1]$. Then it is trivial to check that
$\phi$ is an ortho-order automorphism of ${\mathcal E}(H)$, that is,
beside the order preserving property, $\phi$ also satisfies
\[
\phi(I-E)=I-\phi(E) \qquad (E\in {\mathcal E}(H)).
\]
We can apply the result of our paper \cite{MoPa} to
complete the proof in the case when $\dim H=2$.

Suppose now that $\dim H\geq 3$. We know that $\phi$, when restricted
onto the set of all rank-one projections, is a bijective transformation
preserving orthogonality. It is a celebrated result of Uhlhorn
\cite{Uhl} that in
that case we have a unitary or antiunitary operator $U$ on $H$ such that
\[
\phi(P)=UPU^*
\]
holds for every rank-one projection $P$. Since $\phi(\lambda P)=\lambda
\phi(P)$, it follows that
\[
\phi(\lambda P)=U(\lambda P)U^*
\]
for every $\lambda\in [0,1]$. The
operators of the form $\lambda P$ are exactly the so-called weak atoms
in ${\mathcal E}(H)$ and it is known from \cite{BuGu} that every effect
is equal to the supremum of
the set of all weak atoms which are majorized by the effect in question.
As $\phi$ preserves the order, we thus obtain that
\[
\phi(E)=UEU^*
\]
holds for every $E\in {\mathcal E}(H)$. This completes the proof also in the
present case.
\end{proof}

In our second result which follows we formulate a similar statement
concerning the orthomodular lattice ${\mathcal P}(H)$ of all projections
on $H$. We recall that in the language of effects, the elements of
${\mathcal P}(H)$ are the so-called sharp effects. It should be
emphasized that our result holds only in the case when $\dim
H\geq 3$ (see Remark~\ref{R:kov1}).

\begin{theorem}\label{T:kov1}
Suppose that $\dim H\geq 3$.
Let $\phi :{\mathcal P}(H)\rightarrow{\mathcal P}(H)$ be
a bijective map with the property that
\[
P \leq Q \Longleftrightarrow \phi(P)\leq \phi(Q) \qquad
(P,Q \in {\mathcal P}(H))
\]
and suppose that there are states $D$, $D'$ such that $D$ is not
a scalar operator and
\begin{equation}\label{E:felt}
\text{tr}\, (\phi(P)D')=\text{tr}\, (PD) \qquad (P \in {\mathcal
P}(H)).
\end{equation}
Then there exists an either unitary or antiunitary operator U on H such
that $\phi$ is of the form
\[
\phi(P)=UPU^{*} \qquad (P \in {\mathcal P}(H)).
\]
\end{theorem}

\begin{proof}
As the closed subspaces of $H$ and the projections on $H$ can be
obviously identified, it is clear that our transformation $\phi$ gives
rise to a
so-called lattice automorphism of the lattice of all closed subspaces of
$H$. The main result in \cite{FiLo} states that in the case when
$\dim H\geq 3$, every
such transformation is induced by a bijective semilinear map
$A$ of $H$, and $A$ is bounded and either linear or conjugate-linear if
$\dim H=\infty$. (In fact, the finite dimensional part of this result
was not stated in the corresponding theorem in \cite{FiLo} as it is
just a version of the fundamental theorem of projective
geometry.) For
our transformation $\phi$ this means that \begin{equation}\label{E:kov1}
\phi(P_M)=P_{A(M)}
\end{equation}
holds for every closed subspace $M$ of $H$ ($P_M$ denotes the orthogonal
projection onto $M$).
The main point of our argument is to show that $A$ can be chosen to be
an either unitary or antiunitary operator.

We first show that the semilinear map $A$ corresponding to our
transformation $\phi$ is either linear or conjugate-linear also
in the finite dimensional case. First recall the definition of
semilinearity.
This means that $A$ is additive and $A(\lambda
x)=h(\lambda )A$ holds for every $x\in H$ and $\lambda \in \mathbb C$,
where $h:\mathbb C \to \mathbb C$ is a given ring automorphism.
By \eqref{E:kov1}, for any nonzero $x\in H$ we have
\[
\phi \biggl (\frac{x\otimes x}{\| x\|^2}\biggr )=
\frac{Ax \otimes Ax}{\| Ax\|^2}.
\]
Using \eqref{E:felt} we obtain
\begin{equation}\label{E:f2}
\frac{\text{tr}\, (Ax \otimes Ax \cdot D')}{\| Ax\|^2}=
\frac{\text{tr}\, (x\otimes x \cdot D)}{\| x\|^2}
\end{equation}
which can be rewritten as
\begin{equation}\label{E:kov3}
\frac{\| \sqrt{D'}Ax \|^2}{\| Ax \|^2}=
\frac{\| \sqrt{D}x \|^2}{\| x\|^2}
\end{equation}
for every nonzero $x\in H$. Fix a vector $y\in H$ such that
$\sqrt{D'}Ay \neq 0$.
Replace $x$ by $x+\lambda y$ in \eqref{E:kov3} where
$x\notin \mathbb Cy$ and $\lambda\in \mathbb C$.
Since $A$ is semilinear, we obtain that
\begin{equation}\label{E:f3}
\frac{\Arrowvert\sqrt{D'}(Ax+h(\lambda)Ay) \Arrowvert^2}
{\Arrowvert Ax+h(\lambda)Ay \Arrowvert^2}=
\frac{\Arrowvert\sqrt{D}(x+\lambda y) \Arrowvert^2}
{\Arrowvert x+\lambda y\Arrowvert^2}.
\end{equation}
We now prove that $h$ is bounded in a neighbourhood of $0$. In fact, if
this is not the case, then
there exists a sequence $(\lambda _n)$ in $\mathbb C$ such that
$\lambda _n \rightarrow 0$ and $|h(\lambda _n)| \rightarrow \infty $ as
$n \rightarrow \infty$. Putting these $\lambda_n$'s in the place of
$\lambda$ in \eqref{E:f3} and then taking limit, it follows
that
\[ \frac{\|\sqrt{D'}Ay\|^2}{\| Ay\|^2}=
\frac{\|\sqrt{D}x\|^2}{\| x\|^2}.
\]
Since $y$ was fixed, the left hand side of this equation is constant.
Hence, there is a positive real number $c>0$ such that
\[
\Arrowvert\sqrt{D}x \Arrowvert^2=c\Arrowvert
x\Arrowvert^2
\]
holds for every $x\notin \mathbb Cy$. By continuity we have this
equality for every $x\in H$. This gives us that
\[
\langle Dx,x \rangle =c \langle x,x \rangle  \qquad (x\in H)
\]
which yields $D=cI$. But this case was excluded in the statement and
thus
we obtain that $h$ is bounded in a neighbourhood of 0. It is known
that even such a weak regularity property implies that $h$ is either the
identity or the conjugation on
$\mathbb C$ (one can find a detailed study of the ring automorphisms of
$\mathbb C$ in the book \cite{Kuc} including the mentioned result).
Therefore, we get that $A$ is either linear or conjugate-linear
does not matter what the dimension of $H$ is.

We next prove that $A$ is a scalar multiple of a unitary or antiunitary
operator.
It follows from \eqref{E:f2} that
\[
\frac{\langle D' Ax,Ax \rangle}{\Arrowvert Ax\Arrowvert^2}  =
\frac{\langle D x,x \rangle }{\Arrowvert x\Arrowvert^2}
\]
holds for every $0\neq x\in H$. This implies that
\begin{equation}\label{E:kov4}
\langle D' Ax,Ax \rangle \langle x,x\rangle =
\langle D x,x \rangle \langle Ax,Ax \rangle
\qquad (x\in H).
\end{equation}
Now, fixing $x,y\in H$ for a moment and
replacing $x$ in \eqref{E:kov4} by $\lambda x+y$ where $\lambda$ runs
through $\mathbb{C}$, we
get polynomials of $\lambda$ and $\overline{\lambda}$ on both sides of
\eqref{E:kov4} which are equal
to each other. Therefore, the coefficients of the different
powers of $\lambda$ and $\overline{\lambda}$ in those polynomials are
also equal to each other.
Comparing the coefficients of $\lambda ^2$ we
have
\begin{eqnarray}
\langle D'Ax,Ay \rangle \langle x,y\rangle =
\langle D x,y \rangle \langle Ax,Ay \rangle
\nonumber
\end{eqnarray}
if $A$ is linear and
\begin{eqnarray}
\langle D'Ay,Ax \rangle \langle x,y\rangle =
\langle Dx,y \rangle \langle Ay,Ax \rangle
\nonumber
\end{eqnarray}
if $A$ is conjugate-linear.
Consequently, in both cases we have the implication
\begin{equation}\label{E:kov}
\langle x,y\rangle =0\;\;\;\Longrightarrow \;\;\;
\langle Dx,y \rangle \langle A^*Ax,y\rangle =0.
\end{equation}

As $D$ is compact, according to the spectral theorem of compact
self-adjoint operators, we can write $D=\sum \lambda _i P_i$, where
the $\lambda _i$'s are nonnegative real numbers and the $P_i$'s are
rank-one projections for which there exists an orthonormal basis
$\{e_i \}$ in $H$ such that $P_i e_j=\delta_{ij} e_j$.
Any $x\in H$ can be written as a sum
$x=\sum_i c_i e_i$ with some coefficients $c_i \in \mathbb{C}$.

Suppose that $0\neq x\in H$ is not an eigenvector of $D$. Then
there exist $j,k$ such that
$c_j \neq 0$, $c_k \neq 0$ and $\lambda _j \neq \lambda _k$.
Set $z = \overline{c_k}e_j - \overline{c_j}e_k$. Then we have
\begin{eqnarray}
\langle x,z\rangle=
\langle \sum_i c_i e_i,\overline{c_k}e_j - \overline{c_j}e_k\rangle=
c_j c_k - c_k c_j =0
\nonumber
\end{eqnarray}
and
\begin{eqnarray}
\langle Dx,z\rangle=
\langle x,Dz\rangle=
\langle \sum_i c_i e_i,\lambda _j\overline{c_k}e_j -
\lambda_k\overline{c_j}e_k\rangle=
\lambda _j c_j c_k - \lambda _k c_k c_j \neq 0.
\nonumber
\end{eqnarray}
These imply that for any $y\in H$ with
$\langle x,y\rangle=0$, $\langle Dx,y\rangle=0$ we have
\[
\langle x,y+tz \rangle =0 \quad \text{and} \quad
\langle Dx,y+tz\rangle \neq 0
\]
whenever $t$ is a nonzero real number.
Now, it follows from \eqref{E:kov} that for such $y$ we get
\begin{eqnarray}
\langle A^*Ax,y+tz \rangle=0 \qquad (t\neq 0).
\nonumber
\end{eqnarray}
Taking the limit $t\to 0$, we infer
\[
\langle A^*Ax,y\rangle =0.
\]
By \eqref{E:kov}, we have this equality also in the case when
$\langle x,y\rangle=0$ and $\langle Dx,y\rangle \neq 0$.
Therefore, we have proved
that $\langle A^*Ax,y\rangle =0$ holds for every
$y\in H$ which is orthogonal to $x$. This readily implies
that $A^*Ax$ is a scalar multiple of $x$. Remember that $x$ was an
arbitrary nonzero vector in $H$ which is not an eigenvector of $D$. As
the set of
such vectors is dense in $H$ (this follows from the fact that $D$ is not
scalar), by continuity we deduce that $A^*Ax$ is a scalar multiple of
$x$ for every $x\in H$. This means that, so to say, the linear
operators $A^*A$ and $I$ are locally linearly dependent.

Now, it is a folklore result whose proof requires only elementary linear
algebra that for linear operators of rank at least 2, local linear
dependence implies (global) linear dependence. Hence, we deduce
that $A^*A=\mu I$ holds for some
positive number $\mu $. Therefore, denoting $U=\frac{1}{\sqrt{\mu}}A$,
we have an either
unitary or antiunitary operator $U$ which clearly induces the same
lattice
automorphism on the collection of all closed subspaces of $H$ as $A$.
It is easy to see that this implies
\[
\phi(P)=UPU^* \qquad (P\in {\mathcal P}(H))
\]
and the proof is complete.
\end{proof}

\begin{remark}\label{R:kov1}
One might be interested whether the condition that $D$ is not a scalar
operator can be omitted. We show that this can not be done if $\dim
H\geq 2$. First we note that if $D$ is scalar (which can occur only in
finite dimension), say $D=\lambda I$, then
it is not hard to see that the assumptions in our statement imply that
$D'=D$. In fact, the order preserving property of
the transformation $\phi$ implies that it
preserves the rank-one projections (more generally, the rank-$n$
projections) and then, considering
the equation \eqref{E:felt} for such projections, we obtain that
\[
\text{tr}\, (\phi(P)D')=\text{tr}\, (PD)=\lambda \text{tr}\, P=
\lambda \text{tr}\, \phi(P)=\text{tr}\, (\phi(P) (\lambda I)).
\]
As $\phi(P)$ runs through the set of all rank-one projections,
we deduce from this equality that $D'=\lambda I=D$.
By the rank preserving property of $\phi$, it is now obvious that for
such states (i.e., when $D'=D=\lambda I$) the equality \eqref{E:felt}
does not
represent a proper condition, so all we know then is that $\phi$ is an
order preserving bijection of ${\mathcal P}(H)$. Since,
according to \cite{FiLo}, such a
transformation can be induced by any invertible
semilinear operator, we see that in the treated case
the conclusion in our theorem is no longer true.

We are in the same situation if we omit the condition $\dim
H\geq 3$, that is, our statement does not remain valid
when $\dim H=2$. This is because on such a Hilbert space
the only nontrivial projections are the rank-one projections, so in
that case the order preserving property has no real content, it only
implies that $\phi(0)=0$ and $\phi(I)=I$. Otherwise, on the rank-one
projections $\phi$ can behave arbitrarily.
It is true that the equation \eqref{E:felt} gives some
restriction, but it is just an easy task to find such a
transformation $\phi$ which fulfills the conditions but which does
not preserve the orthogonality between the rank-one projections and
hence can not be written in the form that appears in our theorem.
\end{remark}

% Bibliography
\bibliographystyle{amsplain}

\end{document}